\documentclass[10´pt,a4paper]{article}                     
\usepackage{algorithm2e}
\usepackage[latin1]{inputenc}
\usepackage[T1]{fontenc}
\usepackage{bezier}
\usepackage{dsfont}
\usepackage{mathptmx}      
\usepackage{amsfonts}
\usepackage{amsmath}
\usepackage{amssymb}
\usepackage{amsbsy}
\usepackage{amscd}
\usepackage{amsopn}
\usepackage{amssymb}
\usepackage{amstext}
\usepackage{amsxtra}
\usepackage{graphicx,psfrag}
\usepackage{latexsym}
\usepackage{multirow}
\newtheorem{thm}{Theorem}
\newtheorem{lem}{Lemma}
\newtheorem{wdef}{Definition}

\def\squareforqed{\hbox{\rlap{$\sqcap$}$\sqcup$}}
\def\qed{\ifmmode\else\unskip\quad\fi\squareforqed}
\def\smartqed{\def\qed{\ifmmode\squareforqed\else{\unskip\nobreak\hfil
\penalty50\hskip1em\null\nobreak\hfil\squareforqed
\parfillskip=0pt\finalhyphendemerits=0\endgraf}\fi}}
\smartqed

\DeclareMathAlphabet{\mathcal}{OMS}{cmsy}{m}{n}

\newcommand{\wabs}[1]{\left|#1\right|}

\newcommand{\wcal}[1]{\mathcal{#1}}

\newcommand{\wcn}[1]{{\mathds C}^{#1}}

\newcommand{\wdfc}[2]{{#1}'\!\left(#2\right)}

\newcommand{\wfc}[2]{{#1}\!\left(#2\right)}

\newcommand{\wfloor}[1]{\lfloor {{#1}} \rfloor }

\newcommand{\wgradf}[2]{{\nabla \! #1}\left(#2\right)}
\newcommand{\wgradfc}[2]{{\nabla \! #1}\!\left(#2\right)}

\newcommand{\whessf}[2]{{\nabla^2 \! #1 \left( #2 \right)}}
\newcommand{\whessfc}[2]{{\nabla^2 \! #1 \! \left( #2 \right)}}

\newcommand{\wi}[1]{\wrm{i}}
\newcommand{\wim}[1]{\wfc{\wrm{Im}}{#1}}

\newcommand{\wlc}[1]{\wfc{\wrm{LC}^2}{\wrn{#1}}}

\newcommand{\wlr}[1]{\left( #1 \right)}

\newcommand{\wmod}[2]{ {{#1} \bmod{#2}} }

\newcommand{\wn}{\mathds N}

\newcommand{\wnorm}[1]{\left\| #1 \right\|}

\newcommand{\wre}[1]{\wfc{\wrm{Re}}{#1}}
\newcommand{\wref}[1]{(\ref{#1})}
\newcommand{\wrone}{\mathds R}
\newcommand{\wrn}[1]{{\mathds R}^{#1}}
\newcommand{\wrm}[1]{\mathrm{#1}}

\newcommand{\wseq}[2]{{\left\{ {#1}_{#2}, \ {#2} \in \wn \right\}}}

\newcommand{\wset}[1]{{\left\{ #1 \right\}}}

\newcommand{\wvone}[1]{\mathds{1}_{#1} }

\begin{document}

\title{The divergence of the BFGS and Gauss Newton Methods}
\author{Walter F. Mascarenhas\thanks{
Instituto de Matem\'{a}tica e Estat\'{i}stica, Universidade de S\~{a}o Paulo, 
           Cidade Universit\'{a}ria, Rua do Mat\~{a}o 1010, S\~{a}o Paulo SP, Brazil. CEP 05508-090.
              Tel.: +55-11-3091 5411, Fax: +55-11-3091 6134,  walter.mascarenhas@gmail.com}  
}
\maketitle
\begin{abstract}
We present examples of divergence for the BFGS and Gauss Newton methods.
These examples have objective functions with bounded level sets and
other properties concerning the examples published recently in this journal, like
unit steps and convexity along the search lines. As these other examples,
the iterates, function values and gradients in the new examples
fit into the general formulation in our previous
work {\it On the divergence of line search methods, Comput. Appl. Math. vol.26 no.1 (2007)},
 which also presents an example of divergence for Newton's method.
\end{abstract}

\section{Introduction}
\label{intro}
In the past ten years a few articles have been published, in this journal
and the one mentioned in the abstract, presenting elaborate theoretical examples
of divergence of line search methods. These
methods start from a point $x_0 \in \wrn{n}$ and iterate
according to
\begin{equation}
\label{line_search}
x_{k+1} := x_k + \alpha_k d_k,
\end{equation}
with search directions $d_k \in \wrn{n}$ and parameters $\alpha_k \in \wrone{}$
chosen with the intention that $x_k$
converge to a local minimizer of a function $f: \wrn{n} \mapsto \wrone{}$.
Textbooks  \cite{BERTSEKAS,NOCEDAL} present popular choices for the directions $d_k$
and the parameters $\alpha_k$ and explain why they work in usual circumstances.
The articles \cite{DAIA,DAIB,MEA,MEB} analyze $d_k$'s given by BFGS and conclude that
 these methods may not succeed in extreme situations.
Our article \cite{MEB} presents a similar result for the $d_k$'s corresponding to Newton's method
for minimization and hints that the techniques it describes could be also
applied to other methods, but it does not elaborate on these possible extensions.
Finally, our article \cite{MEC} also offers examples of unexpected behavior for Newton's method.

These examples have subtle points, but none of the them is really perfect.  For instance,
the objective functions in our examples in \cite{MEB} are not explicit and have only
Lipschitz continuous second derivatives. The objective function in \cite{DAIB}
is an explicit polynomial, but it has no local minimizers, its degree is high
and its coefficients are not simple.

The examples are concerned with
the behavior of line search methods in situations that are not
considered in the hypothesis found in textbooks,
or even in the majority of research papers. Usually, textbooks and articles
impose reasonable hypothesis, which are frequently met in practice.
Their spirit is similar to the following theorem
\footnote{For a proof of this theorem, look at equation \wref{mmtsk} in the appendix.}:
\begin{thm}
\label{thm_spd} Consider $n \times n$ matrices $\wseq{M}{k}$, positive
numbers
$\wseq{\alpha}{k}$ and a function
$f: \wrn{n} \mapsto \wrone{}$  with continuous first order derivatives and
bounded level sets.
Suppose the matrices $M_k M_k^t$ are non singular and the iterates $x_k \in \wrn{n}$
are defined by \wref{line_search} with
\begin{equation}
\label{dk_spd}
d_k := - \wlr{M_k M_k^t}^{-1} \wgradfc{f}{x_k}
\end{equation}
and satisfy the first Wolfe condition. If there exists $\overline{\alpha} > 0$ such that
$\alpha_k \geq \overline{\alpha}$ for all $k$
and the matrices $M_k$ are bounded then
$\lim_{k \rightarrow \infty} \wgradf{f}{x_k} = 0$ for every starting point $x_0$.
\end{thm}

The level sets of $f$ are of the form $ \wset{ x \in \wrn{n} \ \wrm{with} \ \wfc{f}{x} \leq a}$
and the first Wolfe condition is the requirement that there exists $\sigma \in (0,1)$ such that
\begin{equation}
\label{first_wolfe}
\wfc{f}{x_{k+1}} \leq \wfc{f}{x_k} + \sigma \wgradfc{f}{x_k}^t \wlr{x_{k+1} - x_k}
\end{equation}
for all $k$. Equation \wref{dk_spd} describes important nonlinear
programming methods, such as steepest descent, BFGS and Gauss Newton.
It also applies to an adaptation of Newton's method in which we take the steepest descent direction
when we detect that the Hessian $\whessf{f}{x_k}$ is not positive definite or
it is almost singular and we would need a very small $\alpha_k$ in order
to obtain a step $s_k = \alpha_k d_k = - \alpha_k \whessfc{f}{x_k}^{-1} \wgradfc{f}{x_k}$ of reasonable size.

The bounded level sets hypothesis in Theorem \ref{thm_spd} implies that
the sequence $x_k$ has a converging subsequence for every $x_0$.
As a result, if the level sets are bounded and we enforce \wref{first_wolfe}
then for every $x_0$ we either are fortunate and have a subsequence converging
to  $x_\infty$  with $\wgradfc{f}{x_\infty} = 0$ or we suffer from one of these two pathologies:
\begin{itemize}
\item[(i)] The parameters $\alpha_k$ in \wref{line_search} get too small, i.e., a subsequence $\alpha_{n_k}$
converges to $0$.
\item[(ii)] The matrices $M_k$ in \wref{dk_spd} are unbounded.
\end{itemize}
As a consequence, for methods like Gauss Newton in which the matrices $M_k$ are
continuous functions of $x_k$, we can only have divergence if the parameters
$\alpha_k$ have a subsequence converging to $0$. Theorem \ref{thm_spd} also
explains why the $\alpha_k$'s in \cite{MEC} converge to $0$.

There are other hypothesis that guarantee the convergence
of line search methods. Some  require
that the search directions do not get almost orthogonal to the gradient or
that Hessians are well conditioned. Others ask for an analytic objective functions
or apply to more elaborate classes of objective functions, as
in \cite{ABSIL,KURDYKA}. Textbooks rely on a combination of
these hypothesis to prove the convergence of the methods they are
concerned with, as we did in Theorem \ref{thm_spd}. They sacrifice
generality for a cleaner exposition of the most common
situations, and are quite right in doing so.
We go in the opposite direction: we explore the
consequences of violating the usual conditions. In theory, we conclude that
methods like BFGS and Gauss Newton may fail if the parameters $\alpha_k$ get
too small or matrices like the $M_k$ above get too large.
This theoretical conclusion disregards rounding errors and the precautions taken in practice.
In fact, items (i) and (ii) above show that the examples will not work in practice.
We cannot handle an unbounded sequence of matrices or arbitrarily small
$\alpha_k > 0$ on a real computer. We would also have trouble working with
the inverses of such matrices.

More than presenting particular examples, we expose the neat analytical,
algebraic and geometric concepts underlying them. There is a subtle relation between the
examples presented here and some techniques to find closed form solutions of
nonlinear differential equations:
they can be both explained in terms of symmetry groups.
By adding a simple term to a nonlinear differential equation
with closed form solutions we may destroy its symmetries and turn it into an equation for which one can
prove that there are no convenient closed form solutions. Examples of
divergence are similar. We did build examples for Newton's method,
BFGS and Gauss Newton using our tools. However,
each example relies on specific features of the method it considers.

This article has four more
sections and one appendix. Section \ref{sec_overview} overviews
the previous examples. It shows how they fit in the framework in \cite{MEB}.
Section \ref{sec_essence} explains the analytical and algebraic underpinnings of the
examples. Section \ref{sec_gauss_newton}  builds an example of divergence for the Gauss Newton method.
Section \ref{sec_bfgs} is about the divergence of BFGS. The appendix
contains proofs and corroborates our claim that
the examples in \cite{DAIB,MEA,MEB} have similar iterates,
function values and gradients. The supplementary material aims to facilitate the reader in using
the software Mathematica to verify the algebra in the examples.

\section{Overview of the examples in \cite{DAIB,MEA,MEB}}
\label{sec_overview}
The examples in \cite{DAIB,MEA,MEB} are based on classical mathematical ideas.
From the numerical point of view,  Powell's work \cite{POWELL} already presents
an interesting analysis of divergence in the same context of functions
with second order derivatives and bounded level sets that we consider here.
From a broader perspective, our ideas and Powell's
are just a natural extension of the mathematical
techniques used to analyze periodic orbits in celestial mechanics in the late 1800s.

The basic ideas behind our work and Powell's are present in the first volume of
Poincar\'{e}'s masterpiece \cite{POINCARE}, which was published in 1892.
In that book Poincar\'{e} uses power series to analyze the convergence of
solutions of differential equations to periodic orbits; Powell uses a
similar technique to analyze the convergence of the iterates of a version
of the conjugate gradient method to a limit cycle. In 1901,
Hadarmard \cite{HADAMARD} looked at the same problem from
a perspective that is quite similar to the one we present in this article.
He was then followed by Cotton \cite{COTTON} and Perron \cite{PERRON}.

\begin{figure}[h]
\includegraphics[bb= 0 0 100 100, scale = 0.1, viewport=0 260 400 500, width=8.5cm, height=2.8cm]{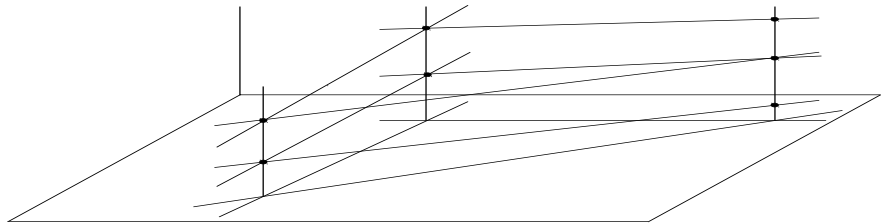}
  \caption{The geometry of divergence: the iterates (black dots) converge to
the vertices of a polygon. The search lines connect consecutive iterates. They converge
to the horizontal lines containing the sides of the polygon.}
\label{fig_divergence}
\end{figure}

From Hadamard's perspective, the dynamics of all the examples of divergence in \cite{DAIB,MEA,MEB} is described
by Figure \ref{fig_divergence}.
In this figure at each step the vertical coordinate is contracted by a factor
$\lambda \in (0,1)$ and the horizontal coordinates are rotated.
The iterates converge to a cycle in the horizontal subspace, which contains no critical points.
In the examples the vertical and horizontal subspaces may have higher dimensions,
the limiting polygon has more vertices and the rotations are replaced by
orthogonal transformations,
but the qualitative picture is the same. \\[0.01cm]

The examples discussed here have features summarized in the following table:
{\small
\begin{center}
\tabcolsep=0.10cm
\begin{tabular}{l|l|l|l }
  \hline
  Feature & \multicolumn{3}{|c}{Article and year of publication } \\[0.01cm]
  \cline{2-4}
               & \cite{DAIB}  -- 2012 & \cite{MEB} -- 2007 & This article -- 2014 \\[0.02cm]
 \hline
Example for Newton's method   & No           & Yes       & No \\
Example for BFGS   & Yes          & Yes       & Yes \\
Example for Gauss Newton      & No           & No        & Yes \\
Bounded level sets            & No           & Yes       & Yes \\[0.04cm]
\hline
Smoothness of the             & Explicit     & Lipschitz continuous & Lipschitz continuous \\
Objective function            & Polynomial   & second derivatives   & second derivatives \\[0.04cm]
\hline
step size ($\alpha_k$)        & One          & Any for Newton, different & One for BFGS\\
                              &              & from one for BFGS         & Converges to $0$  \\
                              &              &                           &  for Gauss Newton \\[0.04cm]
\hline
Convexity along search lines  & Yes & Yes & Yes \\
Armijo and Wolfe conditions   & Yes & Yes & Yes \\
Goldstein condition           & Yes & No & Yes \\
Exact line searches           & Yes & Yes & Yes \\
  \hline
\end{tabular}
\end{center}
}
\noindent
The examples have the same kind of
iterates $x_k$,  function values $f_k = \wfc{f}{x_k}$ and gradients $g_k = \wfc{g}{x_k}$,
which can be written as
\begin{equation}
\label{xk}
x_k = Q^k \wfc{D}{\lambda}^{k}  \overline{x}_k, \hspace{1cm}
f_k = \lambda^{k d_n} \overline{f}_k \hspace{1cm} \wrm{and} \hspace{1cm}
g_k = \lambda^{k d_n} Q^k \wfc{D}{\lambda}^{-k} \overline{g}_k.
\end{equation}
The matrix $Q$ in \wref{xk} is orthogonal and $Q^p = I$ for
a period $p \in \wn{}$. The parameters $\overline{x}_k \in \wrn{n}$,
$\overline{f}_k \in \wrone{}$ and
$\overline{g}_k \in \wrn{n}$
also have period $p$, in the sense that $\overline{x}_{k+p} = \overline{x}_k$,
$\overline{f}_{k+p} = \overline{f}_{k}$ and $\overline{g}_{k+p} = \overline{g}_{k}$.
The matrices $\wfc{D}{\lambda}$ are diagonal. They commute with $Q$ and their diagonal
entries are powers of the parameter $\lambda \in (0,1)$. The constant $d_n$ is equal to
the biggest exponent of $\lambda$ in the diagonal of $\wfc{D}{\lambda}$
(The appendix explains  how the equations in \cite{DAIB,MEA,MEB} fit into
\wref{xk}.)

The parameters in \wref{xk} mix well with Hessians of the form
\begin{equation}
\label{hk}
h_k  = \lambda^{k d_n} Q^k \wfc{D}{\lambda}^{-k} \overline{h}_k \wfc{D}{\lambda}^{-k} Q^{-k}
\end{equation}
and, for the BFGS method, with Hessian approximations of the form
\begin{equation}
\label{bfgs_bk}
B_k = - \sum_{i = 0}^{n-1} \frac{\alpha_{k+i}}{g_{k+i}^t s_{k+i}} g_{k + i} g_{k+i}^t,
\end{equation}
where the $\alpha_k$ are the parameters in \wref{line_search} and also satisfy $\alpha_{k+p} = \alpha_k$.

To build an example of divergence, we express the formulae that define the method we are concerned with,
the Armijo, Goldstein and Wolfe conditions and equation \wref{xk} as a system of equations and inequalities
in $D$, $Q$, $\lambda$, $\overline{x}_k$, $\overline{f}_k$
$\overline{g}_k$, and $B_k$.
We then solve this system of equations and inequalities and interpolate an appropriate objective function at the
$x_k$. Due to the periodicity of $\alpha_k$, $\overline{x}_k$,
$\overline{f}_k$, $\overline{g}_k$ and $\overline{h}_k$ this system consists of a finite number of
equations and inequalities. It is important to realize that we do not need to solve them in closed form.
We can use interval arithmetic and the following version of Moore's Theorem  \cite{RALLY}
to prove that the equations can be solved and get accurate estimates of their solutions
\begin{lem}
\label{lem_moore}
Consider $\overline{x} \in \wrn{n}$, $r > 0$ and
$D = \wset{x \in \wrn{n} \ \wrm{with} \ \wnorm{x - \overline{x}}_1 < r}$.
If $f: D \mapsto \wrn{n}$ has continuous first
derivatives and  the $n \times n$ matrix $A$ and
$a > 0$ are such that
\[
\sup_{1 \leq i \leq n, x \in D} \wnorm{A^t \wgradfc{f_i}{x} - e_i}_1 \wnorm{\wfc{f}{\overline{x}}}_\infty \leq a < 1
\]
and
\begin{equation}
\label{eq_moore}
b := \sup_{1 \leq i \leq n} \wnorm{A^t e_i}_1 \wnorm{\wfc{f}{\overline{x}}}_\infty < r (1 - a)
\end{equation}
then there exists $x^* \in D$ with $\wnorm{x^* - \overline{x}}_\infty \leq b/(1 - a)$ such that
$\wfc{f}{x^*} = 0$.
\end{lem}
\noindent
The proof of Lemma \ref{lem_moore} starts at equation \wref{iKantoC} in the appendix.
Once we obtain estimates for the solution of the equations we can use interval arithmetic to verify the inequalities,
as exemplified in the supplementary material.
Note that all we need to use Lemma \ref{lem_moore} is a good preconditioner $A$ for
the Jacobians of $f$ in $D$. We do not need to estimate Lipschitz constants
for these Jacobians as we would if we were to apply Kantorovich's Theorem \cite{ORTEGA}.

The interpolation processes in \cite{DAIB} and \cite{MEB} are quite different.
In \cite{MEB} we interpolate by extending cubic splines defined along the
search lines to the whole space via Whitney's Extension Theorem, obtaining
an objective function with Lipschitz continuous second derivatives.
The article \cite{DAIB} uses polynomial interpolation. Since the requirement
of Lipschitz continuity of the second derivative is weaker than polynomiality,
our objective functions are more flexible and we can enforce the fundamental condition
of bounded level sets for them. By choosing polynomial
interpolation, \cite{DAIB} is constrained by the lack of flexibility
of analytic functions and, as a consequence, its objective function does not
have local minimizers.

The choice of $d_n = 1$ and appropriate $D$, $Q$, $\lambda$,
$\overline{x}_k$, $\overline{f}_k$, $\overline{g}_k$ and $\overline{h}_k$
in \wref{xk} are basically
all we need to build an example of divergence for Newton's method for minimization
with any
constant positive $\alpha_k$'s. Unfortunately, things are more complicated in the BFGS
method, because we must also handle the matrices $B_k$. Therefore, the
merits of the examples of divergence for the BFGS method in \cite{DAIB} and \cite{MEB}
go beyond their common use of formula \wref{xk} and
the interpolation processes mentioned above. In \cite{DAIB} you will find the end result of
skillful and hard work. In the next sections we present examples based on the
framework developed in \cite{MEB}. None of these examples is trivial.
On the contrary, they are steps towards the noble goal of excellence.

\section{The essence: geometry, algebra and analysis}
\label{sec_essence}
This section outlines how we can build examples of divergence
by combining Whitney's Extension Theorem with the algebra of matrices, if we
are guided by the geometry of Figure \ref{fig_divergence}.
We summarize previous  results so that the reader can have a self contained
view of this process.
The methodical construction of an example of divergence involves two tasks:
\begin{itemize}
\item[(a)] Choosing convenient iterates $x_k$, function values
$f_k$, gradients $g_k$ and Hessians $h_k$ compatible with the method we are concerned with.
\item[(b)] Finding an objective function compatible with the $x_k$, $f_k$, $g_k$ and $h_k$ above.
\end{itemize}
Whitney's Extension Theorem \cite{FEFFERMANA,FEFFERMANB,WHITNEY} is the key ingredient to
impose conditions in the $x_k$, $f_k$, $g_k$ and $h_k$ in item (a) so that we can
perform the interpolation step (b). It is our opinion that
this deep theorem exposes as no other the relation between the nature of functions
with Lipschitz continuous derivatives, the conditions by Armijo, Goldstein
and Wolfe and the theoretical limitations of line searches in spaces of
high dimension. We need the following definition to use Whitney's Extension Theorem
for building examples of divergence:
\begin{wdef}
\label{def_lc2}
We define $\wlc{n}$ as the class of functions $f: \wrn{n} \mapsto \wrone{}$
with Lipschitz continuous second derivatives for which there exists constants
$C$ and $R$, which depend on $f$, such that if
$\wnorm{x} \geq R$
then $\whessf{f}{x}$ is positive definite and
$\wnorm{\whessf{f}{x}^{-1}} \leq C$.
\end{wdef}
We also consider the space $\mathbb{H}_n$ of $n \times n$ symmetric matrices.
Using these concepts we can state the following corollary of Whitney's Extension Theorem
\footnote{Theorem \ref{thm_whitney} is Lemma 6 in page 150 of \cite{MEB}.}:
\begin{thm}
\label{thm_whitney}
Let $E$ be a bounded subset of $\wrn{n}$ and consider functions
$f: E \mapsto \wrone{}$, $g: E \mapsto \wrn{n}$ and
$h: E \mapsto \mathbb{H}_n$. If there exists a constant $M \in \wrone{}$ such that
\begin{eqnarray}
\label{whitney_a}
\wnorm{\wfc{h}{x} - \wfc{h}{y}} & \leq &  M \wnorm{x - y}, \\
\wnorm{\wfc{g}{y} - \wfc{g}{x} - \wfc{h}{x}\wlr{y - x}} &  \leq  &  M \wnorm{x - y}^2, \\
\label{whitney_c}
\wnorm{\wfc{f}{y} - \wfc{f}{x} - \wfc{g}{x}^t \wlr{y - x} - \frac{1}{2} \wlr{y - x}^t \wfc{h}{x}\wlr{y - x}} & \leq &  M \wnorm{x - y}^3,
\end{eqnarray}
then there exists $F \in \wlc{n}$ with
$F(x) = f(x)$, $\wgradfc{F}{x} = \wfc{g}{x}$ and $\whessf{F}{x} = \wfc{h}{x}$
for $x \in E$.
\end{thm}
We could state similar theorems for higher order derivatives, but the algebra needed to express
and handle the consistency conditions analogous to \wref{whitney_a} -- \wref{whitney_c},
and the inequalities they lead to proves to be too complicated.
We believe that the work required to build examples with higher order derivatives would not justify the insights
they would bring. This is why the objective functions
in our examples have only Lipschitz continuous second order derivatives.

Whitney's Extension Theorem exposes a fundamental difference
between our objective functions and the polynomials (or analytic functions) in
examples like \cite{DAIB}.
Analytic functions do not have the extension property described in Theorem \ref{thm_whitney}. They are rigid and cannot be
modified locally. Analytic functions satisfy {\L}ojasiewicz's inequality
\cite{LOJASIEWICZ}, as do the more general families of functions described in
\cite{KURDYKA}. On the other hand, functions in $\wlc{n}$ have bounded level sets and are easy to work with.
By targeting objective functions in this class we do not need to worry about large
$x$ when building examples, because with a little work we can
modify a function with Lipschitz continuous second derivatives to turn
it into an element of $\wlc{n}$. Doing
the same for an analytic function would be a very delicate process,
if feasible at all. This is a fundamental reason why we prefer
functions in $\wlc{n}$ instead of polynomials or analytic
functions. Such choice is also justified as it is
common to find, in textbooks and research papers, theorems
in which the hypothesis asks for a function with Lipschitz
continuous second derivatives and bounded level sets.

We now present Theorem \ref{thm_main}. It is a powerful tool for constructing
examples of divergence for line search methods when combined with matrices
$\wfc{D}{\lambda}$ and $Q$ of the form
\begin{equation}
\label{exampleA}
\wfc{D}{\lambda} =
\left(
\begin{array}{ccc}
 I_a &        0    & 0 \\
0    & \lambda I_b & 0  \\
0    & 0           & \lambda^{d_n} I_c
\end{array}
\right)
\hspace{1cm}
\wrm{and}
\hspace{1cm}
Q =
\left(
\begin{array}{ccc}
Q_a  &   0  & 0   \\
0    & Q_b  & 0   \\
0    &  0   & Q_c \\
\end{array}
\right),
\end{equation}
where $I_i$ stands for the $i \times i$ identity matrix, $a$, $b$ and $c$ are positive integers
and $Q_i$ represents a $i \times i$ orthogonal matrix such that $Q_i^p = I_i$ for some common period $p \in \wn{}$.
There are good reasons for
considering matrices $Q$ and $D$ with three blocks instead of the two blocks used
in some examples in \cite{DAIA,DAIB,MEA,MEB}, and also for considering $a \geq 3$,
$b \geq 2$ and $d_n > 2$ in \wref{exampleA}. These conditions help to ensure
the correct behavior of the limit search lines
\begin{equation}
\label{limit_lines}
\wcal{L}_k  = \wset{\wfc{D}{0} Q^k \wlr{\overline{x}_k + \alpha \overline{s}_k}, \ \ \alpha \in \wrone{}},
\end{equation}
where
\begin{equation}
\label{sk}
\overline{s}_k = Q \wfc{D}{\lambda} \overline{x}_{k+1} - \overline{x}_k.
\end{equation}
With $a \geq 3$ it is quite unlikely that non consecutive limit search lines will cross.
The choice $b \geq 2$ helps to control the rate at which the search lines approach
the limit hyperplane.
Finally, the choice $d_n> 2$ takes care of technical issues regarding the differentiability
of the resulting objective function so that we can satisfy the hypothesis of
Theorem \ref{thm_whitney}.
We can then state the theorem underlying the results in this article
\footnote{This theorem is proved in the last paragraph of the appendix.}:
\begin{thm}
\label{thm_main}
Consider $\lambda \in (0,1)$, the matrices $Q$ and $\wfc{D}{\lambda}$ in \wref{exampleA} and
sequences $\overline{g}_k, \overline{x}_k \in \wrn{n}$, $\overline{f}_k \in \wrone{}$,
$\overline{h}_k \in \mathbb{H}_n$. Suppose that
for a period $p \in \wn{}$ we have $Q^p = I_n$,
$\overline{g}_{k + p} = \overline{g}_k$,
$\overline{x}_{k + p} = \overline{x}_k$, $\overline{f}_{k+p} = \overline{f}_k$ and
$\overline{h}_{k+p} = \overline{h}_k$. If $d_n > 2$ and, for all $k$,
\begin{eqnarray}
\label{distinctLines}
\wfc{D}{0} \overline{s}_k \ \ \wrm{and} \ \ \wfc{D}{0} Q \overline{s}_{k+1}  \ \ \wrm{are \ linearly \ independent}, \\
\label{horzLI}
\wrm{if} \ \wmod{j - k}{p} \not \in \wset{-1,0,1} \ \ \wrm{then} \ \ \wcal{L}_j \cap \wcal{L}_k = \emptyset, \\
\label{vertLI}
\wdfc{D}{0} \overline{x}_k \ \wrm{and} \ \wdfc{D}{0} \overline{s}_k \ \ \wrm{are \ linearly \ independent}, \\
\label{hkk}
\wlr{\overline{h}_k}_{ij} = 0  \ \ \wrm{for} \ d_i + d_j > d_n,
\end{eqnarray}
then there exists $k_0 \in \wn{}$ and $f \in \wlc{n}$  such that for $k \geq k_0$
and $f_k$, $x_k$, $g_k$ and $h_k$ in \wref{xk} and \wref{hk} we have
$\wfc{f}{x_k} = f_k$, $\wgradfc{g}{x_k} = g_k$, $\whessf{f}{x_k} = h_k$.
If we also assume that
\begin{eqnarray}
\label{convexA}
\overline{s}_k^t \overline{g}_k < \lambda^{d_n} \overline{f}_{k+1} - \overline{f}_k <
\lambda^{d_n} \overline{s}_k^t \wfc{D}{\lambda}^{-1} Q \overline{g}_{k+1}, \\
\label{convexB}
 \overline{s}_k^t \overline{h}_k \overline{s}_k > 0 \hspace{0.2cm} \wrm{and} \hspace{0.2cm}
 \overline{s}_k^t Q \wfc{D}{\lambda}^{-1} \overline{h}_{k+1} Q^t \wfc{D}{\lambda}^{-1} \overline{s}_k > 0,
\end{eqnarray}
then $f$ can be chosen to be strictly convex along the search lines, i.e., $s_k^t \whessf{f}{x_k + \alpha s_k} s_k > 0$
for all $\alpha \in \wrone{}$.
\end{thm}

In words, Theorem \ref{thm_main} says that if $d_n > 2$ and the technical conditions \wref{distinctLines}--\wref{convexB} are satisfied
then there exists a function $f$ with bounded level sets and Lipschitz continuous second derivatives which interpolates
the $f_k$, $x_k$, $g_k$ and $h_k$ in \wref{xk}--\wref{hk} for $k$ large enough and is strictly convex along
the search lines. Theorem \ref{thm_main} simplifies the process of building examples of divergence,
because it spares us from the construction of an explicit objective function.
It allows us to concentrate in finding $x_k$, $f_k$, $g_k$ and $h_k$ compatible with our method.
Once we find them we only need to make sure they satisfy the technical conditions \wref{distinctLines}--\wref{convexB}.
As a consequence,  we can explore the iterates in higher dimensions and
observe phenomena which do not occur in lower dimensions.
We can go beyond the two dimensions in which Powell proved the
convergence of the BFGS method and analyze its behavior in situations
in which our fallible lower dimensional intuition may mislead us.
In this exploration we can, and should, take advantage of the modern symbolic
and numerical tools at our disposal, as exemplified in the
supplementary material.  The perception that we can use Theorem
\ref{thm_main} with help of a software allows us to focus on
the creative part in the construction of divergence examples:
the choice of appropriate forms for $\lambda$, $\overline{x}_k$,
$\overline{f}_k$, $\overline{g}_k$ and $\overline{h}_k$ in \wref{xk} and \wref{hk}
and $Q$ in \wref{exampleA}. This choice is the result of our understanding of the method
we are considering and the need to balance the freedom secured by moving to
more dimensions with the complexity of the resulting algebraic problem.

The next sections apply Theorem \ref{thm_main} to build examples for
 Gauss Newton and BFGS. In order to do that we must
look at the particular details that define these methods. We would need
to do the same to every method for which we would like to
apply Theorem \ref{thm_main} to build an example of divergence.
The next three subsections describe technical facts that hold for several methods.
They show that enforcing
conditions like Armijo's, Goldstein and Wolfe and convexity along the
search lines is relatively easy once we get iterates, function values
and gradients that satisfy \wref{xk} and \wref{hk}. The reader
may prefer to skip these technical details and proceed to the next sections.

\subsection{Defining the normalized iterates $\overline{x}_k$ in terms of the normalized steps $\overline{s}_k$}
Usually it is more convenient
to work with the steps $s_k = x_{k+1} - x_k$, and their normalized version
$\overline{s}_k = Q^k \wfc{D}{\lambda}^k \overline{x}_{k+1} - \overline{x}_k$,
instead of the iterates $x_k$ and their normalization $\overline{x}_k$.
Equation \wref{sk} shows how to obtain the $\overline{s}_k$'s from the $\overline{x}_k$'s.
However, we must be cautious when expressing the $\overline{x}_k$'s in
terms of the $\overline{s}_k$'s, since equation \wref{sk} gives rise to a singularity.
In fact, by multiplying $\wref{sk}$ by $Q^{j} \wfc{D}{\lambda}^{j}$ and
recalling that $Q^p = I$ and $\overline{x}_{k+p} = \overline{x}_k$ we obtain
\begin{equation}
\label{x_singular}
\sum_{j = 0}^{p-1} Q^j \wfc{D}{\lambda}^j \overline{s}_{k+j} = \wlr{\wfc{D}{\lambda}^p - I} \overline{x}_k.
\end{equation}
If we decompose $D$ and $Q$ as in \wref{exampleA}, with the corresponding decomposition
\begin{equation}
\label{decompX}
\overline{x}_k = \wlr{\overline{x}_{a,k},\overline{x}_{b,k},\overline{x}_{c,k}}^t \hspace{0.5cm} \wrm{and} \hspace{0.5cm}
\overline{s}_k = \wlr{\overline{s}_{a,k},\overline{s}_{b,k},\overline{s}_{c,k}}^t,
\end{equation}
then \wref{x_singular} leads to these equations:
\begin{equation}
\label{sk_a}
\sum_{j = 0}^{p-1} Q_a^j \overline{s}_{a,k + j} = 0,
\end{equation}
\begin{equation}
\label{xk_bc}
\overline{x}_{b,k} = \frac{1}{\lambda^p - 1} \sum_{j=0}^{p-1} \lambda^j Q_b^j \overline{s}_{b,k + j} \hspace{0.8cm} \wrm{and} \hspace{0.8cm}
\overline{x}_{c,k} = \frac{1}{\lambda^{p d_n} - 1} \sum_{j=0}^{p-1} \lambda^{j d_n} Q_c^j \overline{s}_{c,k + j}.
\end{equation}
We cannot derive $\overline{x}_{a,k}$ from \wref{sk_a}, because the examples are invariant under
translations in the $\overline{x}_{a,k}$'s. However, if \wref{sk_a} holds for $k = 0$ then it holds for all $k$.
Once we enforce \wref{sk_a} for $k = 0$ we can define $\overline{x}_{a,k}$ by
\begin{equation}
\label{xk_a}
\overline{x}_{a,0} = 0 \hspace{0.5cm} \wrm{and} \hspace{0.5cm} k > 0 \hspace{0.1cm} \Rightarrow \hspace{0.1cm}
\overline{x}_{a,k} = Q_a^{-k} \sum_{j = 0}^{k-1} Q_a^j \overline{s}_{a,j}.
\end{equation}
Using induction we can then derive the $\overline{x}_k$ from the $\overline{s}_k$ and prove the following lemma:
\begin{lem}
\label{lem_xk} If, for $k \in \wn{}$, the $\overline{x}_k,\overline{s}_k \in \wrn{n}$ are decomposed as in \wref{decompX}
and satisfy \wref{sk_a} then the $\overline{x}_k$ in
\wref{xk_bc} and \wref{xk_a} are compatible with $\overline{s}_k$ defined in \wref{sk}.
\end{lem}

Using \wref{xk_a} we can write the projection $\wcal{P}_k$ of the limit search line $\wcal{L}_k$
in \wref{limit_lines} in the subspace corresponding to $\overline{x}_{a,k}$ as
\[
\wcal{P}_k := \wset{ \sum_{j=0}^{k-1} Q_a^j \overline{s}_{a,j} + \alpha Q_a^k \overline{s}_{a,k},  \ \ \alpha \in \wrone{}},
\]
under the usual convention that $\sum_{j = 0}^{-1} Q_a^j \overline{s}_{a,k} = 0$.
To verify the hypothesis \wref{horzLI} in Theorem \ref{thm_main} it suffices to show that
$\wcal{P}_{k} \cap \wcal{P}_{k + m} = \emptyset$ for $0 \leq k < p - 1$ and $1 < m < p - 1$. This is
equivalent to saying that there exists no $\alpha,\beta \in \wrone{}$ such that
\[
\sum_{j=0}^{k-1} Q_a^j \overline{s}_{a,j} + \alpha Q_a^k \overline{s}_{a,k} =
\sum_{j=0}^{k + m-1} Q_a^j \overline{s}_{a,j} + \beta Q_a^{k+m} \overline{s}_{a,k + m}
\]
or, equivalently, that there exists no $\alpha,\beta \in \wrone{}$ such that
\begin{equation}
\label{no_intersection}
 \alpha \overline{s}_{a,k} =
\sum_{j=0}^{m-1} Q_a^{j} \overline{s}_{a,k + j} + \beta Q_a^{m} \overline{s}_{a,k + m}.
\end{equation}
In resume, we have proved the following Lemma
\begin{lem}
\label{lem_no_intersection}
If the normalized steps $\overline{s}_k$ are such that for every $0 \leq k < p$ and
$1 < m < p - 1$ there exists no $\alpha$ and $\beta$ satisfying equation \wref{no_intersection}
then the corresponding normalized iterates $\overline{x}_k$ satisfy the hypothesis
\wref{horzLI} of Theorem \ref{thm_main}.
\end{lem}

If $a \geq 3$ the hypothesis of Lemma \ref{lem_no_intersection} will be
satisfied unless the vectors in equation \wref{no_intersection} align in some
unexpected way. As a consequence, we do not need to worry
about this condition as we explore the parameters
that define our example at first. We only need to check
\wref{no_intersection} after we find them.
If, by any chance, the hypothesis of
Lemma \ref{lem_no_intersection} is not satisfied in
the first try, then we should adjust the parameters
slightly, so that this hypothesis holds. If we cannot find a
suitable modification then maybe it would be advisable to
consider whether the method we are considering converges.

\subsection{Convexity along the search lines}
\label{subsec_convex}
The conditions \wref{convexA}--\wref{convexB} enforce convexity along the search lines.
In 2007 we brought up this condition in the abstract of \cite{MEB} in order to make sure that
our examples would choose the only local minimizer along the search line. Under this condition,
the simple algebraic condition $s_k^t g_{k+1} = 0$ guarantees that the iterates in our examples
would be generated by methods that choose a global minimizer along the search line as well as
methods that choose the first local minimizer.
In the examples for Gauss Newton and BFGS in the next sections, and
for other methods that do not use the Hessian of the objective function, it is easy to
enforce \wref{hkk} and \wref{convexB} by decomposing $\overline{h}_k$ in \wref{exampleA} in the block diagonal matrix
\begin{equation}
\label{hk_convex}
\overline{h}_k =
\left(
\begin{array}{ccc}
I_a & 0   & 0 \\
0   & I_b & 0 \\
0   & 0   & 0
\end{array}
\right),
\end{equation}
because in this case \wref{convexB} holds as long as $\overline{s}_{a,k}$ or $\overline{s}_{b,k}$ are not zero.
Since \wref{vertLI} implies that $\overline{s}_{b,k} \neq 0$, by assuming \wref{hk_convex} we do not
need to worry about \wref{convexB}.

\subsection{The conditions by Goldstein, Armijo and Wolfe}

The Goldstein condition \cite{NOCEDAL} requires that there exists $c \in (0,1/2)$ such that
\begin{equation}
\label{goldstein_gen}
\wlr{1 - c} s_k^t g_k \leq f_{k+1} - f_k \leq c s_k^t g_k.
\end{equation}
Using \wref{xk} we can reduce it to
\begin{equation}
\label{goldstein}
(1 - c) \overline{s}_k^t \overline{g}_k \leq \lambda^{d_n} \overline{f}_{k+1} - \overline{f}_k \leq c \overline{s}_k^t \overline{g}_k.
\end{equation}
We can enforce the second inequality in \wref{goldstein} by taking $\overline{f}_k = 1$ for all $k$ and choosing a tiny $c > 0$.
In some cases we can enforce the first inequality in \wref{goldstein} by scaling the $\overline{g}_k$ we already have by $\mu$ in the range
\begin{equation}
\label{Goldstein}
\frac{\lambda^{d_n} \overline{f}_{k+1} - \overline{f}_{k}}{\wlr{1 - c} \max_{1 \leq k \leq p} \overline{s}_k^t \overline{g}_k} \leq \mu \leq
\frac{\lambda^{d_n} \overline{f}_{k+1} - \overline{f}_{k}}{c \min_{1 \leq k \leq p} \overline{s}_k^t \overline{g}_k},
\end{equation}
which is not empty as long as
\[
c \leq \frac{\max_{1 \leq k \leq p} \overline{s}_k^t \overline{g}_k}{\max_{1 \leq k \leq p} \overline{s}_k^t \overline{g}_k +
\min_{1 \leq k \leq p} \overline{s}_k^t \overline{g}_k}.
\]
By imposing the first inequality in the
Goldstein condition \wref{goldstein} we also enforce the first condition in \wref{convexA}.
The second condition in \wref{convexA} follows from the exact line search condition ($s_k^t g_{k+1} = 0$)
when $\overline{f}_{k+1} = \overline{f}_k$, and \wref{convexB} follows from the choice of $\overline{h}_k$ in
\wref{hk_convex}. Therefore, imposing convexity along the search line is not more demanding than enforcing the Goldstein condition
in examples based on Theorem \ref{thm_main}. On the other hand, examples using analytic
function must perform extra work to enforce convexity along the search lines,
as can be noticed by considering the difference in complexity among
the examples with and without this condition present in \cite{DAIB}.

The first Wolfe condition \wref{first_wolfe} is sometimes called Armijo's condition, because Armijo proposed a line search
in which we reduce the step $s_k = x_{k+1} - x_k$ until \wref{first_wolfe} is satisfied.
In our examples the first try for a step already satisfies the first Wolfe condition. Therefore,
the Armijo condition and the first Wolfe condition are equivalent as far as we we are concerned.
Using \wref{xk} we can demonstrate that the first Wolfe condition is satisfied for all positive $\sigma$ smaller than
\[
\sigma_0 = \min_{0 \leq k < p}
\left\{\frac{\lambda^{d_n} \overline{f}_{k+1} - \overline{f}_k}{\overline{s}_k^t \overline{g}_k}  \right\}.
\]
Our examples always have $\overline{s}_k^t \overline{g}_k < 0$. Therefore,
in order for them to satisfy the first Wolfe condition and the Armijo condition it suffices that
$\lambda^{d_n} \overline{f}_{k+1} - \overline{f}_k < 0$ for all $k$.
This can be enforced by taking $\overline{f}_{k} = 1$ for all $k$.

The second Wolfe condition requires that for
some $\beta \in (\sigma, 1)$, where $\sigma$ is the parameter in the first Wolfe condition, we have
\begin{equation}
\label{pre_exact}
\wgradf{f}{x_{k+1}}^t s_k \geq \beta \wgradf{f}{x_k}^t s_k.
\end{equation}
In our examples equation \wref{pre_exact} follows from the descent condition
($s_k^t g_k < 0$) and from the use of exact line searches.
Noticing that \wref{xk} yields
 $s_{k}^t g_{k+1} = \overline{s}_k^t Q \wfc{D}{\lambda}^{-1} \overline{g}_{k+1}$,
we enforce exact line searches  by requiring that
\begin{equation}
\label{exact_line_searches}
\overline{s}_k^t Q \wfc{D}{\lambda}^{-1} \overline{g}_{k+1} = 0.
\end{equation}

\section{An example of divergence for the Gauss Newton method}
\label{sec_gauss_newton}
The purpose of the Gauss Newton
method is to minimize $f: \wrn{n} \mapsto \wrone{}$ given by
\[
\wfc{f}{x} = \frac{1}{2} \sum_{j = 1}^m \wfc{r_j}{x}^2,
\]
where each $r_j$ is a function from $\wrn{n}$ to $\wrone{}$ and $m \geq n$.
The method is defined in terms of the Jacobian matrix $J_r$ of the function
$r: \wrn{n} \mapsto \wrn{m}$ given by $\wfc{r}{x} = \wlr{\wfc{r_1}{x},\dots,\wfc{r_m}{x}}^t$,
whose transpose has the gradients of the $r_j$'s as columns:
\[
\wfc{J_r}{x}^t = \left(\wgradfc{r_1}{x} \mid \wgradfc{r_2}{x} \mid \dots \mid \wgradfc{r_m}{x} \right).
\]
Defining $G_k := \wfc{J_r}{x_k}^t$, and assuming that $G_k G_k^t$ is not singular,
the iterates are defined by
\begin{equation}
\label{def_gauss_newton}
x_{k+1} := x_k - \alpha_k \wlr{G_k G_k^t}^{-1} \wgradfc{f}{x_k} =
x_k - \alpha_k \wlr{G_k G_k^t}^{-1} G_k \wfc{r}{x_k},
\end{equation}
for appropriate $\alpha_k > 0$.
By taking $M_k = G_k$ we see that this method is in the format used in
 Theorem \ref{thm_spd}. If $f$ has bounded level sets and the first Wolfe
condition is satisfied then the matrices $G_k$ are bounded and Theorem \ref{thm_spd}
shows that if the parameters
$\alpha_k$ do not get too close to zero then
$\lim_{k \rightarrow \infty} \wgradf{f}{x_k} = 0$  for every starting point $x_0$.
Therefore, to construct an example of divergence for Gauss Newton we must
allow arbitrarily small $\alpha_k$. Once we accept this fact
we can use Theorem \ref{thm_main} to build an example of divergence for Gauss Newton.

We take $m = n = 7$ and use
Theorem \ref{thm_main} to obtain auxiliary functions $\phi_1,\phi_2,\dots,\phi_7$
with which we define
\begin{equation}
\label{grad_rj}
\wfc{r_j}{x} := \sqrt{\kappa + \wfc{\phi_j}{x}}.
\end{equation}
The parameter $\kappa$ ensures that $\kappa + \wfc{\phi_j}{x} \geq 1$ for all $x$ and $j$. The
$\phi_j$ provided by Theorem \ref{thm_main} are bounded below and we take
\[
\kappa := 1 - \min_{1 \leq j \leq 7} \inf_{x \in \wrn{7}} \wfc{\phi_j}{x}.
\]
As a result we obtain an objective function
\begin{equation}
\label{f_gauss_newton}
\wfc{f}{x} = \frac{1}{2} \wlr{7 \kappa + \sum_{j = 1}^7 \wfc{\phi_j}{x} }.
\end{equation}
The $\phi_j$ are built from the same set of iterates $x_k$ and use the
same search lines as $f$. The $\phi_j$ are convex along the
search lines. Therefore, $f$ is also convex along these lines. Moreover,
$f$ belongs to $\wlc{7}$ because the sum of functions in $\wlc{7}$
is also in $\wlc{7}$. This implies that $f$ has bounded level sets.

We close this section explaining how to use Theorem \ref{thm_main}
to find convenient iterates $x_k$, auxiliary functions $\phi_j$ and parameters
$\alpha_k$ so we can prove that the line search find the only local minimizer
along the search line and satisfy the conditions by Armijo, Goldstein and Wolfe.
We divide the work in four subsections. The first one defines the matrix
$Q$ in \wref{exampleA}, the contraction parameter $\lambda$, the exponent $d_n$ and the iterates $x_k$.
The next one explains how to use Theorem \ref{thm_main} to obtain the functions $\phi_j$.
The third subsection defines the parameters $\alpha_k$ so that the iterates $x_k$ are
compatible with Gauss Newton.
The last subsection shows that the emerging objective function
and iterates satisfy the requirements stated in the abstract.

Finally, we emphasize that although our $\alpha_k$ converge
to zero, the convexity of the objective function along the search lines
and the algebraic condition $s_k^t g_{k+1} = 0$ enforced below imply that these would
be the $\alpha_k$ chosen automatically by an algorithm that follows Powell's
suggestion of choosing the first minimizer along the search line or if
we asked for the global minimizer along the search line.

\subsection{Defining $\lambda$, the matrix $Q$ and the iterates}

We take a constant normalized step
\[
\overline{s}_k := \wvone{7} := \wlr{1,1,1,1,1,1,1}^t.
\]
The matrices $Q$ and $\wfc{D}{\lambda}$ in \wref{xk} are chosen as in \wref{exampleA},
with $a = 4$, $b = 2$, $c = 1$ and
\begin{equation}
\label{q_gauss_newton}
Q_a :=
\left(
\begin{array}{cc}
R_{\pi/3} & 0 \\
0          & R_{\pi/6} \\
\end{array}
\right),
\hspace{1cm}
Q_b := R_{\pi/2}
\hspace{1cm}  \wrm{and} \hspace{1cm}
Q_c := \wlr{-1},
\end{equation}
where
\[
R_\theta =
\left(
\begin{array}{cr}
\cos \theta & - \sin \theta \\
\sin \theta & \cos \theta
\end{array}
\right)
\]
is the counterclockwise rotation by $\theta$.
The parameter $\lambda$ and the period $p$ are defined as
\[
d_n := 3, \hspace{1cm} p := 12 \hspace{1cm} \wrm{and} \hspace{1cm}
\lambda := \sqrt[3]{\frac{1}{1 + \sqrt{3}}}.
\]
As the reader can verify, $\sum_{j = 0}^{11} Q_a^j = 0$. Therefore, Lemma \ref{lem_xk} yields
these normalized iterates
\begin{equation}
\label{xk_gauss_newton}
\overline{x}_k =
\left(
\begin{array}{cc}
\overline{x}_{a,k} \\
\overline{x}_{b,k} \\
\overline{x}_{c,k}
\end{array}
\right)
=
\left(
\begin{array}{cc}
Q_a^{-k} \sum_{j = 0}^{k-1} Q_a^j \wvone{4}  \\
\frac{1}{\lambda^{12} - 1} \sum_{j = 0}^{11} \lambda^j R_{\pi/2}^j \wvone{2} \\
\frac{1}{\lambda^{36} - 1} \sum_{j = 0}^{11} \wlr{- \lambda^3}^j
\end{array}
\right)
=
\left(
\begin{array}{cc}
Q_a^{-k} \sum_{j = 0}^{k-1} Q_a^j \wvone{4}  \\
\wlr{\lambda R_{\pi/2} - I}^{-1} \wvone{2} \\
\frac{-1}{1 + \lambda^{3}}
\end{array}
\right).
\end{equation}
The $D$, $Q$, $\overline{s}_k = \wvone{7}$ and the $\overline{x}_k$ in
the previous equation satisfy the hypothesis \wref{distinctLines} and \wref{vertLI}
of Theorem \ref{thm_main}.
  We end this subsection using Lemma \ref{lem_no_intersection} to verify
the hypothesis \wref{horzLI} in Theorem \ref{thm_main}. In the present case,
equation \wref{no_intersection} reduces to
\begin{equation}
\label{gn_alpha_beta}
\alpha \wvone{4} =
\sum_{j=0}^{m-1} Q_a^{j} \wvone{4} + \beta Q_a^{m} \wvone{4} =
\wlr{I - Q_a}^{-1} \wlr{I - Q_a^m} \wvone{4} + \beta Q_a^m \wvone{4}.
\end{equation}
To verify \wref{horzLI} it would be enough to show that
if $0 < m < 12$ and there are $\alpha,\beta \in \wrone{}$ which satisfy \wref{gn_alpha_beta}
then either $m = 1$ or $m = 11$.
 Let us then prove that these are, indeed, the only two possibilities.
Equations \wref{horzLI}, \wref{q_gauss_newton} and \wref{gn_alpha_beta}
imply that $A \wvone{2}  = B \wvone{2} = 0$, where
\begin{eqnarray}
\label{gn_a}
A & = & \alpha \wlr{I_2 - R_\rho }   - \wlr{I_2 - R_\rho ^m}   - \beta \wlr{I_2 - R_\rho} R_\rho^m, \\
\label{gn_b}
B & = & \alpha \wlr{I_2 - R_\rho ^2} - \wlr{I_2 - R_\rho ^{2m}} - \beta \wlr{I_2 - R_\rho ^2} R_\rho ^{2m},
\end{eqnarray}
for  $\rho = \pi / 6$. The matrices $A$ and $B$ are the sum of a multiple of the identity
and an anti-symmetric matrix. As a consequence, the equalities $A \wvone{2} = B \wvone{2} = 0$ imply that
we actually have $A = B = 0$. Multiplying \wref{gn_a} by $I_2 + R_\rho$ we obtain
\begin{equation}
\label{gn_a2}
\alpha \wlr{I_2 - R_\rho^2}   - \wlr{I_2 - R_\rho ^m} \wlr{I_2 + R_\rho}  - \beta \wlr{I_2 - R_\rho^2} R_\rho^m = 0.
\end{equation}
Subtracting this from \wref{gn_b} we get
\[
- \wlr{I_2 - R_\rho^m} \wlr{R_\rho^m - R_\rho} - \beta \wlr{I_2 - R_\rho^2} R_\rho^m \wlr{R_\rho^m - I_2} = 0.
\]
Since $0 < m < 12$ the matrix $R_\rho^m - I_2$ is non singular. Thus,
$\beta \wlr{I_2 - R_\rho^2} R_\rho^m  = R_\rho - R_\rho^m$
and
\begin{equation}
\label{pre_beta}
\beta \wlr{I_2 - R_\rho^2}  = R_\rho^{1-m} - I_2.
\end{equation}
Since  $\rho = \pi/6$, we have that $\cos 2 \rho = 1/2$  and $\sin 2 \rho = \sqrt{3}/2$. Equating the
entries of the matrix in \wref{pre_beta} to zero we get
\begin{equation}
\label{gn_m}
 \frac{\sqrt{3}}{2} \beta =  \wfc{\sin}{\frac{\wlr{m - 1} \pi}{6}} \hspace{1cm} \wrm{and} \hspace{1cm}
 \frac{1}{2} \beta = \wfc{\cos}{\frac{\wlr{m - 1} \pi}{6}} - 1.
 \end{equation}
It follows that
\[
1 = \wfc{\sin}{\frac{\wlr{m - 1} \pi}{6}}^2 + \wfc{\cos}{\frac{\wlr{m - 1} \pi}{6}}^2  = \frac{3}{4} \beta^2 + \wlr{\frac{1}{2}\beta + 1}^2 = \beta^2 + \beta + 1.
\]
Thus, either $\beta = 0$ or $\beta = -1$. Since $0 < m < 12$, in the case $\beta = 0$ equation \wref{gn_m} implies
 that $m = 1$. Similarly, if $\beta = -1$ then equation \wref{gn_m} implies that $m = 11$ and we are done.

\subsection{Defining the auxiliary functions $\phi_j$}
This subsection explains how to use Theorem \ref{thm_main} to obtain the functions
$\phi_j \in \wlc{7}$. The gradients of these function will be expressed in terms
of the vectors  $e_j \in \wrn{7}$, which satisfy
\[
(e_j)_j = 1 \hspace{1cm} \wrm{and} \hspace{1cm} \wlr{e_j}_i = 0 \ \ \wrm{for} \ \ i \neq j.
\]
We build $\phi_j$ such that
\begin{equation}
\label{good_phi}
\wfc{\phi_j}{x_k} = \lambda^{3k} \hspace{1cm} \wrm{and} \hspace{1cm} \wgradfc{\phi_j}{x_k} = -\lambda^{3k} Q^k \wfc{D}{\lambda}^{-k} e_j,
\end{equation}
by applying Theorem \ref{thm_main} with the $\overline{x}_k$, $\lambda$ and $d_n$ above and
\[
\overline{f}_k = 1 \hspace{1cm} \wrm{and} \hspace{1cm} \overline{g}_k = - e_j.
\]

The $\overline{x}_k$ in the previous subsection satisfy the hypothesis
\wref{distinctLines}--\wref{vertLI} of Theorem \ref{thm_main}. As explained in subsection \ref{subsec_convex},
we can satisfy the hypothesis \wref{hkk} and
\wref{convexB} by taking $\overline{h}_k$ as in \wref{hk_convex}.
The only hypothesis left in order to apply Theorem \ref{thm_main} to obtain $\phi_j$ is
equation \wref{convexA}. In the present case it reduces to
\begin{equation}
\label{gn_convex}
-1 < \lambda^3 - 1 < - \lambda^3 \wvone{7}^t Q^t \wfc{D}{\lambda}^{-1} e_j.
\end{equation}
We now verify the second inequality in the line above in the three possible cases:
\begin{itemize}
\item[(a)] If $1 \leq j \leq 4$ then the second inequality in \wref{gn_convex} follows from
the observation that $\wvone{4}^t Q_a = \wlr{1 + \sqrt{3},1 - \sqrt{3}, \sqrt{3} + 1, \sqrt{3} - 1}^t/2$
and
\[
\lambda^3 - 1 < -0.6 < \frac{1}{2} = - \frac{\lambda^3}{2} \max \wset{1 + \sqrt{3}, 1 - \sqrt{3}, 1 + \sqrt{3}, \sqrt{3} - 1}.
\]
\item[(b)] If $5 \leq j \leq 6$ then the second inequality in \wref{gn_convex} holds because
$\lambda^3 - 1 < -0.6 < - \lambda^2$.
\item[(c)] If $j = 7$ then second inequality in \wref{gn_convex} holds because its right hand side equals one.
\end{itemize}
Therefore, Theorem \ref{thm_main} shows that there exist functions $\phi_j \in \wlc{2}$ which are convex along the lines
$x_k + \alpha s_k$ and satisfy \wref{good_phi}.

\subsection{Defining the $\alpha_k$ for equation \wref{line_search}}
Equations \wref{f_gauss_newton} and \wref{good_phi} show that the gradient of $f$ at $x_k$ is
\begin{equation}
\label{gk_gauss_newton}
g_k :=  \wgradf{f}{x_k} = \frac{1}{2} \sum_{j = 1}^7 \wgradfc{\phi_j}{x} = -
\frac{1}{2} \lambda^{3k} Q^k \wfc{D}{\lambda}^{-k} \wvone{7}.
\end{equation}
Equation \wref{grad_rj} implies that
\[
\wgradf{r_j}{x} = \frac{1}{2 \wfc{r_j}{x}} \wgradf{\phi_j}{x}
\]
and the matrix $G_k = \wfc{J_r}{x_k}^t$ satisfies
\[
G_k =  -\frac{\lambda^{3k}}{2 \sqrt{\kappa + \lambda^{3k}}} Q^k \wfc{D}{\lambda}^{-k}.
\]
Since $Q$ and $\wfc{D}{\lambda}$ commute, $\wlr{G_k G_k^t}^{-1} = 4 \lambda^{-6k} \wlr{\kappa + \lambda^{3k}} \wfc{D}{\lambda}^{2k}$.
Therefore, the search direction $d_k$ satisfies
\[
d_k = - \wlr{G_k G_k^t}^{-1} g_k = 2 \lambda^{-3k} \wlr{\kappa + \lambda^{3k}}Q^k \wfc{D}{\lambda}^k \wvone{7}.
\]
Equations \wref{xk}, \wref{sk} and our choice $\overline{s}_k = \wvone{7}$ yield
$s_k = Q^k \wfc{D}{\lambda}^k \wvone{7}$. Therefore, if we take
\[
\alpha_k := \frac{\lambda^{3k}}{2 \wlr{\kappa + \lambda^{3k}}}
\]
then $s_k = \alpha_k d_k$. Therefore, our iterates are compatible with Gauss Newton
with the $\alpha_k$ above.

\subsection{Verifying the line search conditions}

Equations \wref{xk} and \wref{sk} and our choice $\overline{s}_k = \wvone{7}$ yield
$s_k = Q^k \wfc{D}{\lambda}^k \wvone{7}$. Thus, \wref{gk_gauss_newton} yields
\[
s_k^t g_k = - 7 \lambda^{3k} / 2.
\]
Equation \wref{f_gauss_newton} and our choice $\wfc{\phi_j}{x_k} = \lambda^{3k}$ lead to
$f_k := \wfc{f}{x_k} = 7 \wlr{\kappa + \lambda^{3k}}/ 2$. Thus,
\[
f_{k+1} - f_k = \frac{7}{2} \lambda^{3k} \wlr{\lambda^{3} - 1} = \wlr{1 - \lambda^{3}} s_k^t g_k.
\]
Therefore the first Wolfe condition \wref{first_wolfe} is satisfied for $\sigma = 1 - \lambda^{3}\approx 0.6$.
The Goldstein condition is satisfied for $c = \lambda^3 \approx 0.4$.
Finally, the exact line search condition $s_k^t g_{k+1}$ is satisfied because
$\lambda^3 = 1/(1 + \sqrt{3})$ and according to \wref{gk_gauss_newton}
\[
- 2 s_k^t g_{k+1} = \lambda^3 \wvone{7}^t Q^t \wfc{D}{\lambda}^{-1} \wvone{7} =
 \lambda^3 \wvone{4}^t Q_a \wvone{4} + \lambda^2 \wvone{2}^t R_{\pi/2} \wvone{2} - 1 =
 \lambda^3 \wlr{1 + \sqrt{3}} - 1 = 0.
\]
This completes the construction of the example of divergence for Gauss Newton.

\section{The new example of divergence for the BFGS method}
\label{sec_bfgs}
This section presents a new example of divergence for the
BFGS method. The example shows that this method may fail
even under all the conditions described in the table in section \ref{sec_overview}.
In particular, our examples have bounded level sets, a property
whose far reaching consequences are illustrated in Theorem \ref{thm_spd}.
Time will tell whether it is possible to build an example similar to ours in which
the objective function is a polynomial.

We analyze the BFGS method with exact line searches. In this case
$s_k^t g_{k+1} = 0$ and the BFGS iterates are given by
\begin{equation}
\label{bfgs}
s_k = x_{k+1} - x_k = -\alpha_k B_k^{-1} g_k,
\end{equation}
where $g_k = \wgradf{f}{x_k}$. The positive definite matrices $B_k$ evolve according to
\begin{equation}
\label{bfgs_next}
B_{k+1} := B_k + \frac{\alpha_k}{s_k^t g_k} g_k g_k^t -
\frac{1}{s_k^t g_k} \wlr{g_{k+1} - g_k}\wlr{g_{k+1} - g_k}^t.
\end{equation}
It is convenient to work with $B_k$ of the form \wref{bfgs_bk}, with the additional
requirements that
\begin{equation}
\label{bfgs_skg}
s_k^t g_k < 0 \hspace{1.0cm} \wrm{and} \hspace{1.0cm} s_k^t g_{k+j} = 0 \hspace{1cm} \wrm{for} \hspace{1cm} 1 \leq j < n.
\end{equation}
Equations \wref{bfgs_bk} and \wref{bfgs_skg}
show that $B_k s_k = - \alpha_k g_k$. Therefore, we do not need to worry about \wref{bfgs}.
However we must make sure that the $B_k$ in \wref{bfgs_bk} satisfies
\wref{bfgs_next}. We can achieve this goal by imposing yet another set
of conditions:
\begin{equation}
\label{bfgs_gkpn}
g_{k+n} = \rho_k \wlr{g_{k+1} - g_k},
\end{equation}
for $\rho_k \in \wrone{}$ such that
\begin{equation}
\label{bfgs_akpn}
\alpha_{k+n} = \frac{s_{k+n}^t g_{k+n}}{s_k^t g_k \rho_k^2}.
\end{equation}
Assuming  \wref{bfgs_gkpn} and \wref{bfgs_akpn}, we can use
induction to verify that the matrices \wref{bfgs_bk} satisfy \wref{bfgs_next}. In fact,
if $B_k$ is given by \wref{bfgs_bk} then, using \wref{bfgs_gkpn} and \wref{bfgs_akpn} we obtain
\[
B_{k+1} = - \sum_{i = 0}^{n-1} \frac{\alpha_{k+i}}{g_{k+i}^t s_{k+i}} g_{k + i} g_{k+i}^t +
 \frac{\alpha_k}{s_k^t g_k} g_k g_k^t -
\frac{1}{s_k^t g_k} \wlr{g_{k+1} - g_k}\wlr{g_{k+1} - g_k}^t
\]
\[
=
- \sum_{i = 1}^{n-1} \frac{\alpha_{k+i}}{g_{k+i}^t s_{k+i}} g_{k + i} g_{k+i}^t -
\frac{1}{\rho_k^2 s_k^t g_k} g_{k+n} g_{k+n}^t = - \sum_{i = 1}^{n} \frac{\alpha_{k+i}}{g_{k+i}^t s_{k+i}} g_{k + i} g_{k+i}^t
\]
and $B_{k+1}$ also satisfies \wref{bfgs_skg}.

To build an example of divergence for the BFGS method with
$\alpha_k = 1$ for all $k$ we only need to find $\overline{s}_k$, $\overline{f}_k$, $\overline{g}_k$ and $\rho_k$
which are compatible with the equations \wref{bfgs_skg}--\wref{bfgs_akpn} and the
conditions which allow us to use Theorem \ref{thm_main}. After the experience gained in \cite{MEA,MEB} we
found a better way to parameterize $D$, $Q$, $\overline{x}_k$, $\overline{f}_k$,
$\overline{g}_k$, $\overline{h}_k$ in \wref{xk}--\wref{hk}. Due to a few subtle algebraic
points that we have noticed recently, we now believe it is best to take
\begin{equation}
\label{bfgs_ndp}
n = 9, \hspace{1cm} d_n = 4, \hspace{1cm} p = 16 \times 36 = 576,
\end{equation}
\begin{equation}
\label{example9}
\wfc{D}{\lambda} =
\left(
\begin{array}{ccc}
 I_3 &        0    & 0 \\
0    & \lambda I_2 & 0  \\
0    & 0           & \lambda^4 I_4
\end{array}
\right)
\hspace{1cm}
\wrm{and}
\hspace{1cm}
Q =
\left(
\begin{array}{ccc}
I_3  &   0  & 0   \\
0    & I_2  & 0   \\
0    &  0   & I_4 \\
\end{array}
\right),
\end{equation}
where $I_n$ is the $n \times n$ identity matrix. As we explain below, the choice
of iterates in $\wrn{9}$ and the large period 576 simplifies
the algebra. We define the first 36 normalized iterates, which are then replicated
16 times by symmetry.

Using equation \wref{xk}, parameterized as in
\wref{bfgs_ndp} and \wref{example9}, we reduce the conditions
\wref{bfgs_skg}--\wref{bfgs_akpn} that ensure the compatibility of
our $x_k$ and $B_k$ with the BFGS method with exact line
searches and unity steps to these equations:
\begin{eqnarray}
\label{descent}
\overline{s}_k^t \overline{g}_k & < & 0, \\
\label{rhoConstraint}
\rho_k^2 \overline{s}_k^t \overline{g}_k & = & \lambda^{36} \overline{s}_{k+9}^t \overline{g}_{k+9}, \\
\label{sgOrtho}
\overline{s}_{k}^t Z^j \overline{g}_{k+j} & = & 0 \hspace{0.1cm} \wrm{for} \ j = 1, 2, \dots, 8, \\
\label{gRecursion}
Z^9 \overline{g}_{k+9} & = & \rho_k \wlr{ Z \overline{g}_{k+1} - \overline{g}_k}
\end{eqnarray}
for some $\wset{\rho_k, k \in \wn{}}$ with $\rho_{k+p} = \rho_{p}$ and
\begin{equation}
\label{zLambda}
Z = \wfc{Z}{\lambda} = \lambda^{4} \wfc{D}{\lambda}^{-1} Q =
\left(
\begin{array}{ccc}
\lambda^4 I_3 &    0  & 0 \\
0            &  \lambda^3 I_2 & 0 \\
0            & 0 & I_4
\end{array}
\right).
\end{equation}
A simple inspection of equations \wref{descent}--\wref{gRecursion} proves the following lemma:
\begin{lem}
\label{lemScaleG}
Equations \wref{descent}--\wref{gRecursion} are invariant with respect to
scaling in
$\wset{\overline{g}_k, k \in \wn{}}$,
in the sense that if they are part of a solution of these equations and $\mu > 0$ then
$\wset{\mu \overline{g}_k, k \in \wn{}}$ combined with the same values for the other
parameters also satisfy the same equations.
\end{lem}

The existence of an example of divergence for the BFGS method as claimed in the abstract
is a consequence of the following lemma:

\begin{lem}
\label{lemMain}
There are numbers $\wset{\rho_k, k \in \wn{}}$,
vectors $\wset{\overline{g}_k, k \in \wn{}} \subset \wrn{9}$
and $\wset{\overline{x}_k, k \in \wn{}} \subset \wrn{9}$ such that, for $n$, $p$, $d_n$, $D$ and $Q$ in
\wref{example9}, $\overline{f}_k = 1$, $\overline{h}_k$ in \wref{hk_convex}, $k \in \wn{}$ and
\begin{equation}
\label{pLambda}
\lambda := \sqrt[72]{\frac{1}{1 + \sqrt{2 + \sqrt{2}}}}
\end{equation}
we have $\rho_{k+p} = \rho_k$, $\overline{g}_{k+p} =  \overline{g}_k$ and $\overline{x}_{k+p} = \overline{x}_k$ and
the vectors $\overline{s}_k$ in \wref{sk} satisfy all the conditions
\wref{distinctLines}--\wref{convexB} and \wref{descent}--\wref{gRecursion}.
\end{lem}

We end this section proving Lemma \ref{lemMain}. This demonstration involves the verification
of algebraic identities involving matrices. The supplementary material verifies
them using the software Mathematica.
The idea of the proof is to write the $\overline{g}_k$ and $\overline{s}_k$ as
\begin{equation}
\label{sg}
\overline{g}_{k} := Z^{-k} \Gamma_k e^9_1
 \hspace{0.7cm} \wrm{and} \hspace{0.7cm}
\overline{s}_{k} := -\sigma_k \wlr{\Gamma_{k}^{-1} Z^{k}}^t e^9_1,
\end{equation}
where $Z$ is defined in $\wref{zLambda}$, the $\Gamma_k$ are convenient $9 \times 9$ matrices.
The vector $e^n_i \in \wrn{n}$ has $i$th entry equal to one and the others equal to $0$ and
the $\sigma_k \in \wrn{9}$ are appropriate positive numbers.
We use Lemma \ref{lem_xk} to obtain normalized iterates $\overline{x}_k$ using
\wref{xk} and normalized steps $\overline{s}_k$. The resulting $\overline{x}_k$ are described in equations \wref{xk_bc}
and \wref{xk_a} (We can ignore the matrices $Q_a$, $Q_b$ and $Q_c$ in these
equations because in the present case
they are equal to the identity matrix of the corresponding dimension.)

We consider
\[
u := \lambda^{36} = \sqrt{\frac{1}{1 + \sqrt{2 + \sqrt{2}}}}
\]
and look at $\rho_0,\rho_1, \dots \rho_9$
and $\rho_{18}$ as free parameters. The other $\rho$'s are defined as
\begin{eqnarray}
\rho_k := \rho_{9} \ \ \wrm{for} \ \ k = 10,\dots, 17, & \hspace{0.4cm} &
\rho_k := \rho_{18} \ \ \wrm{for} \ \ k = 19,\dots, 26, \\
\rho_k := \frac{u^2}{\rho_{k - 27} \alpha_1 \alpha_2} \ \ \wrm{for} \ k = 27, \dots, 30,
& \hspace{0.4cm} &
\rho_k := - \frac{u^2}{\rho_{k - 27} \alpha_1 \alpha_2} \ \ \wrm{for} \ k = 31, \dots, 35
\end{eqnarray}
and $\rho_{k} := \rho_{\wlr{\wmod{k}{36}}}$ for $k \geq 36$.
We then define the $9 \times 9$ matrices
\begin{equation}
\label{phik}
\wfc{\Phi}{\rho} := \left(
\begin{array}{cc}
0   & -\rho \\
I_8 & \rho \, e^8_1 \\
\end{array}
\right)
\hspace{1cm} \wrm{and} \hspace{1cm}
\Phi_k := \wfc{\Phi}{\rho_k}.
\end{equation}
These matrices $\Phi_k$ are the KEY part of our arguments, because
the vectors $\overline{g}_k$
satisfy \wref{gRecursion} if and only if the $9 \times 9$  matrices $A_k$ with columns
$Z^{-\wlr{k+i}} \overline{g}_{k + i}$, for $i= 0,\dots 8$, are such that
\[
A_{k+1} = A_k \Phi_k.
\]
Once we grasp how the matrices $A_k$ and $\Phi_k$ are related
it is natural to define
\begin{equation}
\label{psi}
\wfc{\Psi}{\rho} := \prod_{k = 0}^{35} \wfc{\Phi_k}{\rho}
\end{equation}
and search for a vector $\rho$ such that the matrix $\wfc{\Psi}{\rho}^t$ has eigenvalues
$\xi_0, \xi_1, \dots, \xi_8$ given, respectively, by
\begin{equation}
\label{xi}
-u^{4}, \ \ u^{4} e^{ 7 \wi{} \pi/8}, \  \ u^{4} e^{- 7 \wi{} \pi/8}, \  \
 u^3 \wi{}, \ \ - u^3 \wi{}, \ \ e^{\wi{} \pi / 4}, \ \ e^{- \wi{} \pi / 4}, \ \
        -e^{\wi{} \pi / 4} \ \ \wrm{and} \ \  -e^{-\wi{} \pi / 4},
\end{equation}
where $\wi{}$ is the imaginary unit. Once we find $\rho$, we can use the
respective eigenvectors $\nu_0,\nu_1,\dots, \nu_8 \in \wcn{9}$  to
define
\[
\gamma_{2k} := \wre{\nu_{2k}} \in \wrn{9} \hspace{1cm} \ \wrm{and} \hspace{1cm} \gamma_{2k + 1} := \wim{\nu_k} \in \wrn{9}
\]
for $k = 9,\dots,4$.  We then consider the $9 \times 9$ matrix $\Gamma_0$ with rows $\gamma_0,\dots, \gamma_8$
and define
\[
\Gamma_k := \Gamma_0 \prod_{i=0}^{k-1} \Phi_k
\]
(We use  the convention that a product of the form
$\prod_{i = a}^b M_i$ with $b < a$ equals the identity and a sum $\sum_{i = a}^b v_i$ with
$b < a$ is equal to $0$.) Finally we define
\[
\overline{\sigma}_{k} := u^{-\wfloor{k/9}} \wlr{ \prod_{i = 0}^{\wfloor{k/9}-1} \rho_{9 i + \wlr{\wmod{k}{9}}}^2 }.
\]

We now have all the ingredients of equation \wref{sg} and the iterates in
\wref{xk_bc} and \wref{xk_a}. If we look at the eigenvalues $\xi$ in \wref{xi} and
take the $9 \times 9$ block diagonal matrix
\[
\wfc{\Theta}{\lambda} :=
\left(
\begin{array}{ccccc}
- u^{4} &      &                &           & \\
               & \wfc{M}{\xi_1}  &                &           & \\
               &                 & \wfc{M}{\xi_3} &           & \\
               &                 &                & \wfc{M}{\xi_5} & \\
               &                 &                &           & \wfc{M}{\xi_7}
\end{array}
\right)
\]
for
\[
\wfc{M}{r e^{\wi{} \theta }} :=
r \left(
\begin{array}{cc}
\cos \theta & - \sin \theta \\
\sin \theta & \cos \theta
\end{array}
\right)
\]
then we can derive the relation
\begin{equation}
\label{reduced}
\Gamma_0 \wfc{\Psi}{\rho} = \wfc{\Theta}{\lambda} \Gamma_0.
\end{equation}
The identities
\[
\wfc{\Theta}{\lambda} = \wfc{\Theta}{1} Z^{36},
\hspace{1cm}
\sum_{m = 0}^{15} \wfc{\Theta}{1}^m = 0
 \hspace{1cm} \wrm{and} \hspace{1cm}
\prod_{i = 0}^3 \rho_{9 i + k}^2 = u^{4} \ \ \wrm{for} \ k \in \wn{}
\]
and definition \wref{sg} yield
\[
\overline{g}_{36 m + k} = \wfc{\Theta}{1}^m \overline{g}_k
\hspace{1cm} \wrm{and} \hspace{1cm}
\overline{s}_{36 m + k} = \wfc{\Theta}{1}^m \overline{s}_k.
\]
It follows that $\sum_{k= 0}^{p-1} \overline{s}_k^h = 0$ and
we can apply Lemma \ref{lem_xk}. This lemma and equations \wref{xk_bc} and \wref{xk_a}
lead to
\[
\overline{x}_{36 m + k} = \wfc{\Theta}{1}^m \overline{x}_k.
\]
Combining this with $\wfc{\Theta}{1}^{16} = I_9$ and $p = 16 \times 36$ we conclude that
\[
\overline{x}_{k + p} = \overline{x}_k,
\hspace{1cm}
\overline{g}_{k + p} = \overline{g}_k
\hspace{1cm} \wrm{and} \hspace{1cm}
\overline{s}_{k + p} =  \overline{s}_k.
\]
The matrices $\Phi_k$ are such that
\[
\wlr{\prod_{i = 1}^j \Phi_{k+i}} e^9_1 = e^9_{j+1}
\]
for all $k$ and $1 \leq j \leq 8$, because $\Phi_k e^9_i = e^9_{i+1}$ for $1 \leq i \leq 8$.
This implies \wref{sgOrtho} and, similarly, the remaining conditions in
\wref{descent}--\wref{gRecursion} can be verified by plugging \wref{sg} into them.

Therefore, all we need to produce $\overline{g}_k$, $\overline{s}_k$ and $\overline{x}_k$
that would satisfy the compatibility conditions \wref{descent}--\wref{gRecursion}
 is to find $\rho_0,\dots,\rho_9$ and $\rho_{18}$
in such way that the matrix $\wfc{\Psi}{\rho}$ in $\wref{psi}$ has the eigenvalues $\xi$ in $\wref{xi}$.
The Mathematica script in the supplementary material proves that
these $\rho_k$'s exist using Lemma \ref{lem_moore}.
It also finds bounds on them, computes the corresponding eigenvectors,
$\overline{g}_k$, $\overline{s}_k$ and $\overline{x}_k$ using interval arithmetic
and shows that these parameters satisfy the geometric constraints
\wref{distinctLines}--\wref{vertLI}. The proof of Lemma \ref{lemMain} is thus completed.

\appendix
\section{Technicalities}
This appendix begins with an explanation
as to why the form of the iterates, function values and gradients in reference
\cite{DAIB} was already described in \cite{MEB}. We then prove Lemma \ref{lem_moore}
and, finally, the theorems.

Let us then see how several equations in \cite{DAIB} and \cite{MEB} correspond to
\wref{xk}. In \cite{DAIB} and the example for Newton's method in \cite{MEB}
the parameter $d_n$ is equal to $1$.
In the example for the BFGS method in \cite{MEB} it is equal to $3$. Equation \wref{xk} corresponds to
equations (10)--(12) in \cite{MEB}. It is generalized in equations (34)--(37) of \cite{MEB}.
\cite{DAIB} adds a constant $f^*$ to $f$, but this is irrelevant. It calls $\lambda$ by $t$ and thinks
in terms of the steps $\delta_k = x_{k+1} - x_k$, so that $\delta_k = Q^k \wfc{D}{\lambda}^k \overline{\delta}$
in  \cite{DAIB}'s equation (2.1). This is just an affine change of coordinates of \wref{xk}
and only a different notation for the corresponding $s_k$'s used in \cite{MEB}.
The matrix
$M$ in \cite{DAIB}'s equation (2.2) is equal to the matrix $Q \wfc{D}{\lambda}$, whose $k$th power
multiplies $\overline{x}_k$ in equation \wref{xk}. The matrix $P$ in  \cite{DAIB}'s equation
(2.5) is equal to the matrix $\lambda Q \wfc{D}{\lambda}^{-1}$ whose $k$th power multiplies
$\overline{g}_k$ in equation \wref{xk}. Therefore, the iterates, the function values
and the gradients are described in essentially the same way in the examples in \cite{DAIB,MEA,MEB},
in terms of orthogonal matrices $Q$ scaled by powers of $\lambda$ (or $t$) via
the matrices $\wfc{D}{\lambda}$ or their inverses.

{\bf Proof of Lemma \ref{lem_moore}}. Let us define $\delta_0 = 0$ and $x_0 = \overline{x}$
and consider $\delta_k$ and $x_k$
defined inductively by $\delta_k = A \wfc{f}{x_{k-1}}$ and $x_k = x_{k-1} - \delta_k$.
The lemma will be proved if we show that
\begin{equation}
\label{iKantoC}
\wnorm{\delta_{k}}_\infty \leq a^{k-1} b,
\hspace{0.5cm}
 \wnorm{x_{k+1} - \overline{x}} \leq \frac{1 - a^k}{1 - a} b
\hspace{0.5cm} \wrm{and} \hspace{0.5cm}
\wnorm{\wfc{f}{x_k}}_\infty \leq a^k \wnorm{\wfc{f}{\overline{x}}}_\infty,
\end{equation}
because these bounds imply that $x_k$ converges to $x_\infty$ with
$\wnorm{x_\infty - \overline{x}}_\infty \leq b / (1-a)$ and $\wfc{f}{x_\infty} = 0$.
Let us then prove \wref{iKantoC} by induction. Equation \wref{iKantoC}
certainly holds for $k = 0$. Assuming that \wref{iKantoC} holds for $k$, we conclude
from \wref{eq_moore} and definition of $\delta_k$ that
\[
\wnorm{\delta_{k+1}}_\infty \leq \sup_{1 \leq i \leq n} \wnorm{A^t e_i}_1 \wnorm{\wfc{f}{x_k}}_\infty
\leq a^{k} \wnorm{A^t e_i}_1 \wnorm{\wfc{f}{\overline{x}}}_\infty \leq a^k b
\]
and the first bound in \wref{iKantoC} holds for $k + 1$. The second bound on
\wref{iKantoC} follows from the analogous bound for $\wnorm{x_k - \overline{x}}_\infty$ and the
bound in $\delta_{k+1}$ above. It shows that $x_{k+1} \in D$
and the segment $S$ connecting $x_k$ to $x_{k+1}$ is contained in $D$.
As a result, the Mean Value Theorem for the function $f_i: D \mapsto \wrone{}$ implies that,
for some $\xi \in S$,
\[
\wabs{\wfc{f_i}{x_{k+1}}} =
\wabs{\wfc{f_i}{x_k - A \wfc{f}{x_k}}} =
\wabs{\wfc{f_i}{x_k} - \wgradfc{f_i}{\xi}^t A \wfc{f}{x_k}}
\]
\[
=
\wabs{\wlr{A^t \wgradfc{f_i}{\xi} - e_i}^t \wfc{f}{x_k}}
\leq \wabs{A^t \wgradfc{f_i}{\xi} - e_i}_1 \wnorm{\wfc{f}{x_k}}_\infty
\leq a \wnorm{\wfc{f}{x_k}}_\infty.
\]
Thus, $\wabs{\wfc{f}{x_{k+1}}}_\infty \leq a \wnorm{\wfc{f}{x_k}}_\infty
\leq a^{k+1} \wnorm{\wfc{f}{\overline{x}}}_\infty$
and we are done.\qed{}

{\bf Proof of Theorem \ref{thm_spd}.} We start by rewriting \wref{line_search} and
\wref{dk_spd} as
\begin{equation}
\label{mmtsk}
M_k M_{k}^t s_k  = - \alpha_k g_k,
\end{equation}
for $s_k := x_{k+1} - x_k$ and $g_k := \wgradfc{f}{x_k}$.
Equation \wref{mmtsk} shows that $s_k^t g_k = - \wnorm{M_k^t s_k}^2 / \alpha_k \leq 0$
and the first Wolfe condition \wref{first_wolfe} yields
\begin{equation}
\label{proof_df}
f_k := \wfc{f}{x_k} \geq f_{k+1}.
\end{equation}
Since $f$ has bounded level sets this implies that
\[
x_k \in K := \wset{x \in \wrn{n} \ \wrm{with} \ \wfc{f}{x} \leq \wfc{f}{x_0}}.
\]
The set $K$ is compact and the matrices $M_k$ are bounded. Since $f$ has continuous first derivatives
there exists a constant
$\kappa$ such that $\wnorm{g}_k \leq \kappa$,
$\wnorm{M_k} \leq \kappa$, $\wnorm{x_k} \leq \kappa$ and $\wnorm{s_k} \leq \kappa$.

To prove that $\lim_{k\rightarrow \infty} g_k = 0$ we use the well known
result that if a bounded sequence $\wseq{u}{k} \subset \wrone{}$ is such
that all its converging subsequences $u_{n_k}$ converge to zero then
$u_k$ itself converges to $0$.
Let us then consider a subsequence $g_{n_k}$ such that
$\lim_{k \rightarrow \infty} \wnorm{g_{n_k}} = L$ and show that $L = 0$.
Equation \wref{mmtsk} leads to $\wnorm{g_k} \leq \kappa^3 / \alpha_k$.
Thus, if some sub subsequence of $\alpha_{n_{k_j}}$ converges to $+\infty$ then
$L = \lim_{j \rightarrow \infty} \wnorm{g_{n_{k_j}}} = 0$
and we are done. We can then assume that there exists $A$ such that
 $\alpha_{n_k} \leq A$ for all $k$.
Equations \wref{first_wolfe}, \wref{mmtsk} and \wref{proof_df} yield
\[
\wnorm{M_{n_k}^t s_{n_k}}^2 \leq - \alpha_{n_k} s_{n_k}^t g_{n_k} \leq
\alpha_{n_k} \wlr{ \wfc{f}{x_{n_k}} - \wfc{f}{x_{n_k+1}}} / \sigma
\leq \alpha_{n_k} \wlr{ \wfc{f}{x_{n_k}} - \wfc{f}{x_{n_{k+1}}}} / \sigma
\]
\[
\leq A \wlr{ \wfc{f}{x_{n_k}} - \wfc{f}{x_{n_{k+1}}}} / \sigma
\]
and
\[
\sum_{k = 0}^\infty \wnorm{M_{n_k}^t s_{n_k}}^2 \leq \frac{A}{\sigma} \wlr{ \wfc{f}{x_{n_0}} - \inf_{x \in K} \wfc{f}{x} }.
\]
Therefore,  $\lim_{k \rightarrow \infty} M_{n_k}^t s_{n_k} = 0$. If follows that
$g_{n_k} = - M_{n_k} M_{n_k}^t s_{n_k} / \alpha_{n_k}$ converges to $0$, because the $M_{n_k}$ are bounded and $\alpha_{n_k} \geq \overline{\alpha}$.
\qed{}

{\bf Proof of Theorem \ref{thm_main}.} This theorem is an specialization of Theorem 4 in page 145 of \cite{MEB}.
The reader can prove Theorem 3
by using Theorem 4 in \cite{MEB} by realizing that the hypothesis of Theorem 3 corresponds
to a simplified version of the more general concept of seed, which is described in definition
7 in page 142  of \cite{MEB} (There is a typo in the end of the item 3
in this definition. It should read like "...$\wfc{D}{0}\overline{s}_r$  and $\wfc{D}{0}Q \overline{s}_{r+1}.$",
with a $Q$ between $\wfc{D}{0}$ and $\overline{s}_{r+1}$.) The reader will note that our choice for
the diagonal exponents in \wref{exampleA} ($0$ in the first block, $1$ in the second and $d_n > 2$)
leads to vacuous sub itens (c) and (d) in the fouth item in the definition of seed. As a result,
we do not need to worry about these items. Morevoer, we only need to worry about item (b) in
the case $i = j = n$. This is why we have hypothesis \wref{hkk} in Theorem \ref{thm_main}.
The reader can now read the many details and involved proving Theorem \cite{MEB}.
This proof is very technical and there is no point repeating it here.

\end{document}